\documentclass[12pt]{article}

\usepackage{amsmath, epsfig, cite}
\usepackage{amssymb}
\usepackage{amsfonts}
\usepackage{latexsym}
\usepackage{float}
\usepackage{color}
\usepackage{mathrsfs}

\newtheorem{thm}{Theorem}[section]

\newtheorem{lem}[thm]{Lemma}

\newcommand{\pf}{\noindent{\it Proof} }

\numberwithin{equation}{section}

\newcommand{\qed}{{\hfill$\square$}\medskip}

\begin{document}
	%\linenumbers
	\begin{center}
		{\large\bf   Further $q$-Supercongruences from Singh's Quadratic Transformation  }
	\end{center}
	\vskip 2mm \centerline{Wei-Wei Qi}
	
	\begin{center}
		{\footnotesize MOE-LCSM, School of Mathematics and Statistics, Hunan Normal University, Hunan 410081, P.R. China\\[5pt]
			{\tt wwqi2022@foxmail.com} \\[10pt]
		}
	\end{center}
	
	\vskip 0.7cm \noindent{\bf Abstract.} In this paper, we investigate some $q$-congruences  for truncated ${}_{4}\phi_3$ series by using Singh's quadratic transformation and the   “creative microscoping” method (introduced by Victor J. W. Guo and Zudilin in $2019$).

	\vskip 2mm
	\noindent{\it MR Subject Classifications}: 33D15, 11A07, 11B65	
	
	\section{Introduction} 
	
	In $1997$, Van Hamme \cite[(H.2)]{w5} showed the following interesting supercongruence: for any odd prime $p$,	
	\begin{align}
		\begin{aligned}
			\sum_{k=0}^{(p-1)/2}\frac{(\frac{1}{2})^3}{k!^3} \equiv 
			\begin{cases}
				\Gamma_p(1/4)^4 \pmod{p^2}\quad&{if \quad p\equiv 1\pmod 4},\\
				0	 \pmod{p^2} \quad&{if\quad  p\equiv 3\pmod 4},\label{in-1}
			\end{cases} 
		\end{aligned}
	\end{align}		
	where $(a)_n=a(a+1)\dots(a+n-1)$ stands for the Pochhammer symbol and $\Gamma_p(x)$ is the $p$-adic Gamma function. It is clear that $(1/2)/(k!)\equiv 0 \pmod{p}$ for $(p+1)/2 \leq k \leq p-1$, hence, \eqref{in-1} remains valid when the sum is taken over $k$ from 0 to $p-1$. 
	
	Over the past few decades, $q$-analogues of supercongruences constitute an interesting research topic, which  have been widely studied by many scholars. A number of distinct generalizations of \eqref{in-1} have been given in (\cite{w7}, \cite{w8}, \cite{w9} and so on). Particularly, Guo and Zudilin 	\cite[Theorem 2]{w7} showed a $q$-analogue of \eqref{in-1} as follows: modulo $\Phi_n(q)^2$,
	\begin{align}
		\begin{aligned}
			\sum_{k=0}^{(n-1)/2}\frac{(q^2;q^4)_k^2(q^2;q^4)_k}{(q^2;q^2)_k^2(q^4;q^4)_k}q^{2k} \equiv 
			\begin{cases}
				\frac{(q^2;q^4)_{(n-1)/4}^2}{(q^4;q^4)_{(n-1)/4}^2}q^{(n-1)/2} \quad&{if \quad n\equiv 1\pmod 4},\\
				0	 \quad&{if\quad  n\equiv 3\pmod 4}.
			\end{cases} 
		\end{aligned} \label{in-2}
	\end{align}
	Here  and in what follows, the \textit{$q$-shifted factorial} (\textit{$q$-Pochhammer symbol}) is given by 
	\begin{align*}
		(a;q)_0=1, \quad (a; q)_n=(1-a)(1-aq)\dots(1-aq^{n-1}) \quad for \quad n=1,2,\dots
	\end{align*}
	The \textit{$q$-integer} $[n]_q$ is defined by $[n]_q=\frac{1-q^n}{1-q}=(1+q+\dots+q^{n-1})$. In addition, the $n$th cyclotomic polynomial is given by
	\begin{align*}
		\Phi_n(q)=\prod_{\substack{1\le k \le n\\[3pt](n,k)=1}}(q-\zeta^k),
	\end{align*}
	where $\zeta$ denotes a primitive $n$th root of unity.
	
	For more $q$-congruences and $q$-supercongruences derived from transformations for basic hypergeometric series, together with various techniques, one may consult (\cite{w6}, \cite{w6-1}, \cite{w1}, \cite{w2}, \cite{w4} and so on).
	Recently, Guo and Zhao \cite{w1} proposed the following $q$-supercongruence related to the aforementioned \eqref{in-2}: Let $d\geq 2$ be an integer and $x$ an indeterminate. Let $n$ be a positive integer $n\equiv 1 \pmod{2d}$ . Then
	\begin{align}
		\sum_{k=0}^{(n-1)/d}\frac{(q;q^d)_k^2(x;q^d)_k}{(q^d;q^d)_k(q^{d+2};q^{2d})_k}q^{dk}\equiv \sum_{k=0}^{(n-1)/(2d)}\frac{(q;q^{2d})_k^2(x^2;q^{2d})_k}{(q^{2d};q^{2d})_k(q^{d+2};q^{2d})_k}q^{2dk} \pmod{\Phi_n(q)^2}.  \label{in-3}
	\end{align}	
	For $d=2$, \eqref{in-3} reduces to
	\begin{align}
		\sum_{k=0}^{(n-1)/2}\frac{(q;q^2)_k^2(x;q^2)_k}{(q^2;q^2)_k(q^{4};q^{4})_k}q^{2k}\equiv \sum_{k=0}^{(n-1)/4}\frac{(q;q^{4})_k^2(x^2;q^{4})_k}{(q^{4};q^{4})_k^2}q^{4k} \pmod{\Phi_n(q)^2},  \label{in-4}
	\end{align}	
	where $n\equiv 1\pmod{4}$.	
	we observe that $(q;q^2)_k^2$ contains the factor $(1-q^n)^2$ for $  (n-1)/2 \leq k\leq n-1$. Meanwhile, the denominators $(q^2;q^2)_k(q^4;q^4)_k$ of \eqref{in-4} remain coprime to $\Phi_n(q)$. Therefore, it is natural to see that the left-hand side summation in \eqref{in-4} can be extended to $n-1$. Similarly, the right-hand side of \eqref{in-4} can also be extended to  $n-1$. Namely, we have
	\begin{align*}
		\sum_{k=0}^{n-1}\frac{(q;q^2)_k^2(x;q^2)_k}{(q^2;q^2)_k(q^{4};q^{4})_k}q^{2k}\equiv \sum_{k=0}^{n-1}\frac{(q;q^{4})_k^2(x^2;q^{4})_k}{(q^{4};q^{4})_k^2}q^{4k} \pmod{\Phi_n(q)^2}.
	\end{align*}	
	
	Motivated by the above, we further obtain the following consequences:		
	\begin{thm}
		For a positive $d\geq 2$, a positive integer $n$ with  $n \equiv 1 \pmod{2d}$,   and an indeterminate $x$, we have
		\begin{align}
			\sum_{k=0}^{n-1}\frac{(q;q^d)_k^2(x;q^d)_k}{(q^d;q^d)_k(q^{d+2};q^{2d})_k}q^{dk}\equiv \sum_{k=0}^{n-1}\frac{(q;q^{2d})_k^2(x^2;q^{2d})_k}{(q^{2d};q^{2d})_k(q^{d+2};q^{2d})_k}q^{2dk} \pmod{\Phi_n(q)^2}. \label{th-2}
		\end{align}	
		Let $d\geq 3$ be an integer, $n$ a positive integer with  $n \equiv -1 \pmod{2d}$, and $x$  an indeterminate. Then,
		\begin{align}
			\sum_{k=0}^{n-1}\frac{(q^{-1};q^d)_k^2(x;q^d)_k}{(q^d;q^d)_k(q^{d-2};q^{2d})_k}q^{dk}\equiv \sum_{k=0}^{n-1}\frac{(q^{-1};q^{2d})_k^2(x^2;q^{2d})_k}{(q^{2d};q^{2d})_k(q^{d-2};q^{2d})_k}q^{2dk}  \pmod{\Phi_n(q)^2}. \label{th-1}
		\end{align}
		
	\end{thm}			
	
	Taking $x=\pm q^2$ in \eqref{th-2}, we arrive at the following $q$-congruence: for a positive integer $n\equiv 1 \pmod{4}$, we have
	\begin{align}
		\sum_{k=0}^{n-1}\frac{(q;q^2)_k^2}{(q^{4};q^4)_k}q^{2k}\equiv \sum_{k=0}^{n-1}\frac{(q;q^{4})_k^2}{(q^{4};q^{4})_k}q^{4k} \pmod{\Phi_n(q)^2},\label{th-2-0}
	\end{align}
	and
	\begin{align}
		\sum_{k=0}^{n-1}\frac{(q;q^2)_k^2}{(q^{2};q^2)_k^2}q^{2k}\equiv \sum_{k=0}^{n-1}\frac{(q;q^{4})_k^2}{(q^{4};q^{4})_k}q^{4k} \pmod{\Phi_n(q)^2}. \label{th-2-1}
	\end{align}
	Inspired by \eqref{th-2-0} and \eqref{th-2-1}, we deduce the following conclusion.
	\begin{thm}
	(i)	Let $d\geq 2$ be an integer and  $n>1$  a positive odd integer with $n\equiv -1 \pmod d$. Then 
		\begin{align}
			\sum_{k=0}^{\frac{n+1}{d}}\frac{(q^{-1};q^d)_k}{(q^d;q^d)_k}q^{dk}\equiv (-1)^{\frac{n+d+1}{d}}[n]_qq^{\frac{n^2-dn+d-1}{2d}} \pmod{\Phi_n(q)^2}. \label{th-2-2}
		\end{align}		
		For a positive integer $d$ and  a positive odd integer $n>1$ with $n \equiv 1 \pmod{2d}$, we have
		\begin{align}
			\sum_{k=0}^{\frac{n-1}{2d}}\frac{(q;q^{2d})_k}{(q^{2d};q^{2d})_k}q^{2dk}\equiv 0 \pmod{\Phi_n(q)}. \label{th-2-3}
		\end{align}	
	
(ii) Let $d\geq 2$ be an integer, x ($x\neq 0$) an indeterminate, and  $n>1$  a positive odd integer with $n\equiv -1 \pmod d$. Then 
\begin{align}
	\sum_{k=0}^{\frac{n+1}{d}}\frac{(q^{-1};q^d)_k(x;q^d)_k}{(q^d;q^d)_k(q^{d-1}x;q^d)_k}q^{dk}\equiv (-1)^{\frac{n+d+1}{d}}[n]_qq^{-\frac{n^2+dn-d-1}{2d}}\frac{((qx)^{-1};q^d)_{(n+1)/d}}{(x^{-1};q^d)_{(n+1)/d}} \pmod{\Phi_n(q)^2}. \label{th-2-4}
\end{align}		
For a positive integer $d$, an indeterminate x ($x\neq 0$) and  a positive odd integer $n>1$ with $n \equiv 1 \pmod{2d}$, we have
\begin{align}
	\sum_{k=0}^{\frac{n-1}{2d}}\frac{(q;q^{2d})_k(x;q^{2d})_k}{(q^{2d};q^{2d})_k(q^{2d+1}x;q^{2d})_k}q^{2dk}\equiv 0 \pmod{\Phi_n(q)}. \label{th-2-5}
\end{align}	
\end{thm}

	In particular, letting $n$ be a prime and $q\rightarrow 1$ in \eqref{th-2-2} yields the following congruence: for odd prime $p\equiv 3\pmod{4}$, 
	\begin{align*}
		\sum_{k=0}^{\frac{p+1}{4}}\frac{(-\frac{1}{4})_k}{k!}\equiv (-1)^{\frac{p-3}{4}}p \pmod{p^2}. 
	\end{align*}

	\begin{thm}
		Let $d\geq 2$ be an integer, $n$ a positive integer, $x$ an indeterminate and $s\in \{-1,1\}$ with  $n \equiv s \pmod{2d}$. Then 
		\begin{align}
			\sum_{k=0}^{n-1}\frac{(x;q^d)_k(q^s;q^d)_k}{(q^d;q^d)_k}q^{dk}\equiv \sum_{k=0}^{n-1}\frac{(x^2;q^{2d})_k(q^s;q^{2d})_k}{(q^{2d};q^{2d})_k}q^{2dk} \pmod{\Phi_n(q)}. \label{th-3}
		\end{align}		
	\end{thm}
	
	Letting $x=0$, $n=p$ be a prime and $q\rightarrow 1$ in \eqref{th-3}, we immediately obtain the following result: for an integer $d\geq 2$, $s\in \{-1,1\}$ and a prime $p$ with $p\equiv s \pmod{2d}$, then
	\begin{align*}
		\begin{aligned}
			\sum_{k=0}^{p-1}\frac{(\frac{s}{d})_k}{k!}\equiv 	\sum_{k=0}^{p-1}\frac{(\frac{s}{2d})_k}{k!} \pmod{p}.
		\end{aligned}
	\end{align*}
	
	\begin{thm}
		For a positive integer $d$, an odd positive integer $n$, and an indeterminate $x$ and $y$, we have
		\begin{align}
			\sum_{k=0}^{n-1}\frac{(x;q^d)_k(y;q^d)_k}{(xyq^d;q^{2d})_k}q^{dk}\equiv \sum_{k=0}^{n-1}\frac{(x;q^{2d})_k(y;q^{2d})_k}{(xyq^{d};q^{2d})_k}q^{2dk} \pmod{\Phi_n(q)}. \label{th-5}
		\end{align}
	\end{thm}
	
The rest of the paper is organized as follows.  In Section $2$, we  present a $q$-congruence with a parameter $a$. In Section $3$, we prove the main theorems.
	
	\section{A $q$-Congruence With Parameter $a$}
	Recall  that Singh's quadratic  transformation \cite[Appendix (III.21)]{w3} be stated as follows: 	
	\begin{align}
		\begin{aligned}
			{}_{4}\phi_3\left[ \begin{array}{c}
				a^2, b^2, c, d \\
				ab\sqrt{q}, -ab\sqrt{q}, -cd
			\end{array} ; q, q \right]
			={}_{4}\phi_3\left[ \begin{array}{c}
				a^2, b^2, c^2, d^2 \\
				a^2b^2q, -cd, -cdq
			\end{array} ; q^2, q^2 \right], \label{s-1}
		\end{aligned}
	\end{align}	
	where the basic hypergeometric series ${}_{r+1}\phi_r$ (cf. \cite{w2}) is defined as	
	\begin{align*}
		\begin{aligned}
			{}_{r+1}\phi_r\left[ \begin{array}{c}
				a_1, a_2, \dots, a_{r+1} \\
				b_1, b_2, \dots, b_r
			\end{array} ; q, z \right]
			=\sum_{k=0}^{\infty}\frac{(a_1;q)_k(a_2;q)_k\dots(a_{r+1};q)_k}{(q;q)_k(b_1;q)_k\dots(b_r;q)_k} z^k.
		\end{aligned}
	\end{align*}	
	It is clear that the $d=0$ case of \eqref{s-1} reduces to	
	\begin{align}
		\begin{aligned}
			{}_{3}\phi_2\left[ \begin{array}{c}
				a^2, b^2, c \\
				ab\sqrt{q}, -ab\sqrt{q}
			\end{array} ; q, q \right]
			={}_{3}\phi_2\left[ \begin{array}{c}
				a^2, b^2, c^2 \\
				a^2b^2q, 0
			\end{array} ; q^2, q^2 \right]. \label{s-2}
		\end{aligned}
	\end{align}	
	The transformation \eqref{s-2} may be considered as a $q$-analogue of Gauss' quadratic transformation
	\begin{align*}
		{}_{2}F_1(2a,2b;a+b+\frac{1}{2};z)=	{}_{2}F_1(a,b;a+b+\frac{1}{2};4z(1-z)).
	\end{align*}	
	
	Now, we use \eqref{s-2} to derive the following $q$-congruence involving parameters $a$.	
	\begin{lem}
		Let $d\geq 3$ be an integer, $n$  a positive integer, $x$  an indeterminate and $s\in \{-1,1\}$ with  $n \equiv s \pmod{2d}$. Then, modulo $(1-aq^n)(a-q^n)$,
		\begin{align}
			\sum_{k=0}^{n-1}\frac{(aq^s,q^s/a;q^d)_k(x;q^d)_k}{(q^d;q^d)_k(q^{d+2s};q^{2d})_k}q^{dk}\equiv \sum_{k=0}^{n-1}\frac{(aq^s,q^s/a;q^{2d})_k(x^2;q^{2d})_k}{(q^{2d};q^{2d})_k(q^{d+2s};q^{2d})_k}q^{2dk}. \label{s-3}
		\end{align}
		For  a positive integer $n$ with  $n \equiv 1 \pmod{4}$ and an indeterminate $x$. Then, modulo $(1-aq^n)(a-q^n)$,
		\begin{align}
			\sum_{k=0}^{n-1}\frac{(aq,q/a;q^2)_k(x;q^2)_k}{(q^2;q^2)_k(q^{4};q^{4})_k}q^{2k}\equiv \sum_{k=0}^{n-1}\frac{(aq,q/a;q^{4})_k(x^2;q^{4})_k}{(q^{4};q^{4})_k^2}q^{4k}. \label{s-3-1}
		\end{align}
	\end{lem}
	\pf.  Applying the parameters substituting $q\rightarrow q^d$, $a= q^{(s-n)/2}$, $b= q^{(s+n)/2}$ and $c=x$ in Singh's transformation \eqref{s-2}, we get
	\begin{align}
		\begin{aligned}
			{}_{3}\phi_2\left[ \begin{array}{c}
				q^{s-n}, q^{s+n}, x \\
				q^{(d+2s)/2}, -q^{(d+2s)/2}
			\end{array} ; q^d, q^d \right]
			={}_{3}\phi_2\left[ \begin{array}{c}
				q^{s-n}, q^{s+n}, x^2 \\
				q^{d+2s}, 0
			\end{array} ; q^{2d}, q^{2d} \right]. \label{s-4}
		\end{aligned}
	\end{align}	
	Since $(q^{s-n};q^d)_k=0$ for $(n-s)/d < k < (n-1)$ and $(q^{s-n};q^{2d})=0$ for $(n-s)/(2d) < k< (n-1)$. Hence, we have
	\begin{align}
		\sum_{k=0}^{n-1}\frac{(q^{s-n};q^d)_k(q^{s+n};q^d)_k(x;q^d)_k}{(q^d;q^d)_k(q^{d+2s};q^{2d})_k}q^{dk}=\sum_{k=0}^{n-1}\frac{(q^{s-n};q^{2d})_k(q^{s+n};q^{2d})_k(x^2;q^{2d})_k}{(q^{2d};q^{2d})_k(q^{d+2s};q^{2d})_k}q^{2dk}.
	\end{align}	
	It follows that  both sides of \eqref{s-3} are equal for $a=q^n$ and $a=q^{-n}$. Meanwhile $1-aq^n$ and  $a-q^n$ are coprime polynomials in $q$. Therefore, the $q$-congruence \eqref{s-3} is true modulo $(1-aq^n)(a-q^n)$. Similarly, setting $d=2$ and $s=1$ in \eqref{s-4}, we  obtain \eqref{s-3-1}. Thus, we complete the proof. \qed

	\section{Proofs of the Theorems}
	
	\textit{Proof of Theorem $1.1$}. For $d\geq 2$ and $n\equiv s\pmod{2d}$, we know that $gcd(2d,n)=1$. On the one hand, $(q^d; q^d)_k$ and $(q^{2d};q^{2d})_k$ are coprime to $\Phi_n(q)$ for $0\leq k \leq n-1$. On the other hand, when $d\neq 2$, the smallest positive integer $k$ satisfying $(q^{s};q^{2d})_k\equiv 0$, $(q^{s};q^{d})_k\equiv 0$ and $(q^{d+2s};q^{2d})_k\equiv 0$ modulo $\Phi_n(q)$ are $(n-s+2d)/(2d)$, $(n-s+d)/d$ and $((2+d)n+d-2s)/(2d)$, respectively. Noting that
	\begin{align*}
		\frac{n-s+2d}{2d} < \frac{n-s+d}{d} < \frac{(2+d)n+d-2s}{2d}.
	\end{align*}	
	It follows that for $ 0\leq k\leq n-1$ and $d\geq 2$, the polynomials $(q^d;q^d)_k(q^{d+2s};q^{2d})_k$ and $(q^{2d};q^{2d})_k(q^{d+2s};q^{2d})_k$
	are relatively prime to $\Phi_n(q)$. Meanwhile, we have the fact $q^n\equiv 1\pmod{\Phi_n(q)}$. Therefore,  setting $a=1$ in \eqref{s-3} gives \eqref{th-1}. Letting $a=1$ in \eqref{s-3-1} and combining \eqref{s-3} yields \eqref{th-2}.
	
	\textit{Proof of Theorem $1.2$}. Recall that the $q$-Chu-Vandermonde summation \cite[Appendix (II.6)]{w3} can be written as:
	\begin{align}
		\begin{aligned}
			{}_{2}\phi_1\left[ \begin{array}{c}
				a, q^{-n} \\
				c
			\end{array} ; q, q \right]
			=\frac{(c/a;q)_n}{ (c;q)_n}a^n. \label{ss-0-0}
		\end{aligned}
	\end{align}	
	For a positive integer $d\geq 2$. Making the parameter substitutions $q\rightarrow q^d$, $a\rightarrow q^{n-1}$, $n\rightarrow (n+1)/d$, and $c\rightarrow q^{-1}$ in \eqref{ss-0-0}, we get
	\begin{align}
		\begin{aligned}
			\sum_{k=0}^{(n+1)/d}&\frac{(q^{n-1};q^d)_k(q^{-n-1};q^{d})_k}{(q^d;q^d)_k(q^{-1};q^d)_k}q^{dk}\\
			&=q^{(n^2-1)/d}(q^{-n};q^d)_{(n+1)/d}/(q^{-1};q^d)_{(n+1)/d}\\
			&=-(-1)^{(n+1)/d}[n]_qq^{(n^2-dn+d-1)/(2d)}. \label{ss-0-1}
		\end{aligned}
	\end{align}
	Considering that for $d\geq 2$ and $n\equiv -1\pmod{d}$, the denominators $(q^d;q^d)_k$ of \eqref{ss-0-1}
	are relatively prime to $\Phi_n(q)$. Hence, we immediately obtain the desired result \eqref{th-2-2}  from \eqref{ss-0-1}.
	
	Moreover, letting $q\rightarrow q^{2d}$, $a\rightarrow q^{1-n}$, $n\rightarrow (n-1)/(2d)$, and $c\rightarrow q$ in \eqref{ss-0-0}, we have
	\begin{align}
		\begin{aligned}
			\sum_{k=0}^{\frac{n-1}{2d}}\frac{(q^{1-n};q^{2d})_k^2}{(q^{2d};q^{2d})_k(q;q^{2d})_k}q^{2dk}=\frac{(q^n;q^{2d})_{(n-1)/(2d)}}{(q;q^{2d})_{(n-1)/(2d)}}q^{-\frac{(n-1)^2}{2d}}.
			\label{ss-0-2}
		\end{aligned}
	\end{align}
	Since the denominators of \eqref{ss-0-2} are all coprime with $\Phi_n(q)$, and $(q^n;q^{2d})_{(n-1)/(2d)}$  contains a factor $1-q^n$. Therefore, we get \eqref{th-2-3} by utilizing \eqref{ss-0-2}.
	
In view of $q$-Saalach\"{u}tz identity \cite[Appendix(II.12)]{w3}:	
\begin{align}
	\begin{aligned}
		{}_{3}\phi_2\left[ \begin{array}{c}
		q^{-n},	a, b  \\
			c, q^{1-n}ab/c
		\end{array} ; q, q \right]
		=\frac{(c/a;q)_n(c/b;q)_n}{ (c;q)_n(c/(ab);q)_n}. \label{ss-0-3}
	\end{aligned}
\end{align}		
Then, performing the parameter substitutings $q\rightarrow q^{d}$, $a\rightarrow q^{n-1}$, $n\rightarrow (n+1)/{d}$, $c\rightarrow x$, and $c\rightarrow q^{-1}$ in \eqref{ss-0-3}, we have	
\begin{align}
	\begin{aligned}
		\sum_{k=0}^{(n+1)/d}&\frac{(q^{n-1};q^d)_k(q^{-n-1};q^{d})_k(x;q^d)_k}{(q^d;q^d)_k(q^{-1};q^d)_k(q^{d-1}x;q^d)_k}q^{dk}\\
		&=\frac{(q^{-n};q^d)_{(n+1)/d}((qx)^{-1};q^d)_{(n+1)/d}}{(q^{-1};q^d)_{(n+1)/d}(q^{-n}x^{-1};q^d)_{(n+1)/d}}\\
		&=-(-1)^{(n+1)/d}[n]_qq^{-\frac{(n-1)(n+1-d)}{2d}}\frac{((qx)^{-1};q^d)_{(n+1)/d}}{(q^{-n}x^{-1};q^d)_{(n+1)/d}}. \label{ss-0-4}
	\end{aligned}
\end{align}	
In addition, taking $q\rightarrow q^{2d}$, $a\rightarrow q^{1-n}$, $n\rightarrow (n-1)/(2d)$, $b\rightarrow x$, and $c\rightarrow q$ in \eqref{ss-0-3}, we obtain	
	\begin{align}
	\begin{aligned}
		\sum_{k=0}^{\frac{n-1}{2d}}\frac{(q^{1-n};q^{2d})_k^2(x;q^{2d})_kq^{2dk}}{(q^{2d};q^{2d})_k(q;q^{2d})_k(q^{2d-2n+1}x;q^{2d})_k}=\frac{(q^n;q^{2d})_{(n-1)/(2d)}(q/x;q^{2d})_{(n-1)/(2d)}}{(q;q^{2d})_{(n-1)/(2d)}(q^n/x;q^{2d})_{(n-1)/(2d)}}.
		\label{ss-0-5}
	\end{aligned}
\end{align}	
The proofs of \eqref{th-2-4} and \eqref{th-2-5} then follow from the identities \eqref{ss-0-4} and \eqref{ss-0-5}. Since \eqref{th-2-4} and \eqref{th-2-5} can be derived in a similar manner to the proofs of \eqref{th-2-2} and \eqref{th-2-3}, we omit their detailed proofs here. 	
	
	\textit{Proof of Theorem 1.3}.
	Let $d\geq 2$ be an integer, $n$ a positive integer, $x$ a indeterminate and $s\in \{-1,1\}$ with  $n \equiv s \pmod{2d}$. Then, modulo $(1-aq^n)$,
	\begin{align}
		\sum_{k=0}^{n-1}\frac{(x;q^d)_k(aq^s;q^d)_k}{(q^d;q^d)_k}q^{dk}\equiv \sum_{k=0}^{n-1}\frac{(x^2;q^{2d})_k(aq^s;q^{2d})_k}{(q^{2d};q^{2d})_k}q^{2dk}. \label{s-5} 
	\end{align}			
	Making the parameter substitutions $q\rightarrow q^d$, $a=q^{(s-n)/2}$, $b=0$ and $c=x$ in Singh's transformation \eqref{s-2},  for $n\equiv s \pmod{2d}$ we obtain 
	\begin{align*}
		\begin{aligned}
			{}_{2}\phi_1\left[ \begin{array}{c}
				q^{s-n}, x \\
				0
			\end{array} ; q^d, q^d \right]
			={}_{2}\phi_1\left[ \begin{array}{c}
				q^{s-n}, x^2 \\
				0
			\end{array} ; q^{2d}, q^{2d} \right]. 
		\end{aligned}
	\end{align*}	
	Noting that $(q^{s-n};q^d)_k=0$ for $(n-s)/d < k < (n-1)$ and $(q^{s-n};q^{2d})_k=0$ for $(n-s)/(2d) < k< (n-1)$, we see that both sides of \eqref{s-5} are equal for $a=q^{-n}$. Therefore, the $q$-congruence holds modulo $1-aq^n$. Taking $a=1$ in \eqref{s-5}, we obtain \eqref{th-3}.

	Likewise, we have a parametric generalization of Theorem $1.4$.
	\begin{lem}
		For a positive integer $d$, an odd positive integer $n$, an indeterminate $x$, and a non-negative indeterminate $m$. Then, modulo $\Phi_n(q)$
		\begin{align}
			\sum_{k=0}^{n-1}\frac{(m,q^d)_k(x;q^d)_k(y;q^d)_k}{(-mq^d;q^d)_k(xyq^d;q^{2d})_k}q^{dk}\equiv \sum_{k=0}^{n-1}\frac{(m^2;q^{2d})_k(x;q^{2d})_k(y;q^{2d})_k}{(-mq^d;q^d)_{2k}(xyq^{d};q^{2d})_k}q^{2dk}. \label{ss-0}
		\end{align}	
	\end{lem}	
	\pf. Letting $q\rightarrow q^d$, $a\rightarrow \sqrt{x}$, $b\rightarrow \sqrt{y}$, $c\rightarrow q^{d(1-n)}$, and $d\rightarrow m$ in \eqref{s-1}, we have
	\begin{align*}
		\begin{aligned}
			{}_{4}\phi_3\left[ \begin{array}{c}
				x, y, q^{d(1-n)}, m \\
				\sqrt{xyq^d}, -\sqrt{xyq^d}, -mq^{d(1-n)}
			\end{array} ; q^d, q^d \right]
			={}_{4}\phi_3\left[ \begin{array}{c}
				x, y, q^{2d(1-n)}, m^2 \\
				xyq^d, -mq^{d(1-n)}, -mq^{d(1-n)+1}
			\end{array} ; q^{2d}, q^{2d} \right], 
		\end{aligned}
	\end{align*}		
	which can be written as follows
	\begin{align}
		\begin{aligned}
			&\sum_{k=0}^{n-1}\frac{(x;q^d)_k(y;q^d)_k(q^{d(1-n)};q^d)_k(m;q^d)}{(q^d;q^d)_k(-mq^{d(1-n)};q^d)_k(xyq^d;q^{2d})_k}q^{dk}\\
			&=\sum_{k=0}^{n-1}\frac{(x;q^{2d})_k(y;q^{2d})_k(m^2;q^{2d})_k(q^{2d(1-n)};q^{2d})_k}{(q^{2d};q^{2d})_k(-mq^{d(1-n)};q^d)_{2k}(xyq^d;q^{2d})_k}q^{2dk} \label{ss-1}
		\end{aligned}	
	\end{align}	
	Since $q^n\equiv 1 \pmod{\Phi_n(q)}$,  for $0\leq k \leq n-1$ and non-negative $m$, the polynomials $(q^d;q^d)_k(-mq^d;q^d)_k$ and $(q^{2d};q^{2d})_k(-mq^d;q^d)_{2k}$ are relatively prime to $\Phi_n(q)$. Hence, In light of  \eqref{ss-1} we get \eqref{ss-0}, which confirms Lemma $3.1$. \qed
	
	\textit{Proof of Theorem $1.4$}. Letting $m=0$ in \eqref{ss-0}, we obtain the desired $q$-congruence \eqref{th-5}. Now, the proofs of all Theorems are finalized. \qed

\end{document}